\documentclass[11
pt,reqno,fleqn]{amsart}
\usepackage{amsfonts,amsmath,amsthm,amssymb}
\usepackage{pdfsync}

\usepackage{mathrsfs}

\newtheorem{thm}{Theorem}

\newtheorem{lemma}{Lemma}

\newcommand{\R}{\ensuremath{\mathbb{R}}}
\newcommand{\Z}{\ensuremath{\mathbb{Z}}}
\newcommand{\C}{\ensuremath{\mathbb{C}}}

\newcommand{\bfa}{\ensuremath{\mathbf{a}}}
\newcommand{\bfq}{\ensuremath{\mathbf{q}}}

\newcommand{\<}{\ensuremath{\lesssim}}

\newcommand{\eps}{\ensuremath{\varepsilon}}
\newcommand{\la}{\ensuremath{\lambda}}

\newcommand{\ds}{\ensuremath{\widetilde{d\sigma}}}

\newcommand{\eq}{\begin{equation}}
\newcommand{\ee}{\end{equation}}

\newcommand{\p}{\ensuremath{\mathfrak{p}}}

\numberwithin{equation}{section}

\usepackage[paper=a4paper,margin=1in]{geometry}

\newgeometry{top=1in,bottom=1in,right=0.75in,left=.75in}

\begin{document}
\title{A note on discrete spherical averages over sparse sequences}
\author{Brian Cook}
\date{}

\begin{abstract} This note presents an example of an increasing sequence $(\la_l)_{l=1}^\infty$ such that the maximal operators associated to normalized discrete spherical convolution averages \[
\sup_{l\geq 1}\frac{1}{r(\la_l)}\left|\sum_{|x|^2=\la_l}f(y-x)\right|\]
for functions $f:\Z^n\to\C^n$ are bounded on $\ell^p$ for all $p>1$ when the ambient dimension $n$ is at least five. 
\end{abstract}

\thanks{
The author was supported in part by NSF grant DMS1147523}

\maketitle

\section{Introduction}

Let $|x|^2=x_1^2+...+x_n^2$ for a given $n$ and define the convolution averages \[
A_\la f(y)=\frac{1}{r(\la)}\sum_{|x|^2=\la}f(y-x)\]
for functions $f:\Z^n\to\C$. Here $r(\la)=\#\{x\in\Z^n:|x|^2=\la\}$ denotes the number of lattice points on the discrete sphere of radius $\la^{1/2}$ where $\la\geq1$ is an integer. An integral analogue of Stein's spherical maximal theorem \cite{Stein} is given in \cite{MSW}. This result states that the maximal operator\[
S_* f(y)=\sup_{\la\geq 1}\,|A_\la f(y)|\]
is a bounded operator on $\ell^p(\Z^n)$, meaning that there is a constant $C$ which possibly depends on $n$ and $p$ such that \[
\sum_{y\in\Z^n}|S_* f(y)|^p\leq C^p \sum_{y\in\Z^n}| f(y)|^p,\]
if and only if $n\geq 5$ and $p>n/(n-2)$. 

In the original setting, over $\R^n$ that is,  Stein's result shows $L^p(\R^n)$-boundedness when $n>2$ and $p>n(n-1)$ of the maximal operators $\mathcal{A}_*$ defined by  \[
\mathcal{A}_* f=\sup_{\la>0}|f*\sigma_\la|,\]
the supremum being over all real $\la>0$,  where \[
\mathcal{A}_\la f(x)=\int f(y-\la x)d\sigma (x).\]
Here $\sigma$ denotes the normalized Euclidean surface measure on the sphere of radius one in $\mathbb{R}^n$. The case $n=2$ is also a well known result and is due to Bourgain \cite{BouCirc}.  The range of $p$ in the results of Stein and Bourgain may not be enlarged, however this is no longer the case when the supremum is restricted to range only over a sparse sequence of radii. The principle example is for the maximal operators\[
\mathcal{D}_* f=\sup_{j\in\Z}|f*\sigma_{2^j}|,\] operators which are known to be bounded on  $L^p(\R^n)$ for all $p>1$. This may be found in \cite{Calderon} and \cite{Guido}, and a more general result of this type may be found in \cite{Rubio}.

A natural question then concerns  whether or not the same phenomenon is present in the discrete setting as well.  More specifically, for a given sparse sequence $\la_l$ does one have $\ell^p$-boundedness for the maximal operators\eq\label{1.1}
\sup_{l\geq 1}\,|A_{\la_l} f(y)|\ee
for some range $p\in(p_0,\infty)$ with $p_0<n/(n-2)$? This question is first  considered by Hughes in \cite{Kevin} and he provides some results in this direction (and further results can be found in \cite{Lacey}). However an example attributed to  Zienkiewicz shows that attaining a result with $p_0=1$ is impossible in any dimension for discrete spherical averages simply on the condition that the sequence is lacunary. The point of this note is to show that, for some sequence of radii at least, we do recover $\ell^p$-boundedness or the operator \eqref{1.1} in the full range $p>1$. 

Our main observation is that, along certain sequences, spherical averages behave in a manner similar to the averages over homogenous hypersufaces considered by the author in \cite{cook}. A prime example in the homogeneous case concerns conic averages defined by the forms $\p(x)=x_1^2+...+x_{n-1}^2-x_{n}^2$, which are given by \[
C_N f(y)=\frac{1}{N^{n-2}}\sum_{\p(x)=0;\,|x_1|,...,|x_n|\leq N}f(y-x).\] 
The maximal averages $C_*f=\sup_{N\geq0}|C_N f|$ are indeed bounded on $\ell^p$ for any $p>1.$

Let $\la_l=2^l!$ and define the maximal operator \[
A_*f(y)=\sup_{l>0} |A_{\la_l}f(y)|.\]
We show the following.

\bigskip
\begin{thm}\label{thm1}
If $n\geq 5$ and $p>1$ then $||A_*||_{\ell^p\to\ell^p}<\infty$.
\end{thm}
\bigskip

It is worth pointing out some easily verified facts about  our sequence $\la_l$. The first takes the form of a  sparsity condition, namely that $\la_l^2<\la_{l+1}$. This  gives an order of growth which is roughly doubly dyadic, meaning somewhat  comparable to the dyadic subsequence of the dyadic sequence of integers. The second condition is a quantification of the observation that $\la_l$ is (essentially) zero modulo $q$ for all natural numbers $q$. The precise condition used  is that $\la_l\equiv 0 \mod(2^i!)$ for all $l\geq i$. The latter condition, while not necessary as one could easily interlace a finite number of translates of this sequence, is the reason there is no conflict with the examples of Zienkiewicz.

\bigskip
\section{Outline}

The proof follows the argument of \cite{cook}, which is itself a blend of the arguments  of Bourgain in (\cite {Boup>1}, section 7) and Magyar, Stein, and Wainger \cite{MSW}. The current setup is different from that of \cite{cook}, although with the sequence chosen in the introduction the differences end up being minimal and much of that work will apply to our current setup with minimal modifications. 

On the $n$-dimensional $\Pi^n$, identified with the fundamental domain $[0,1]^n$, we begin by defining the Fourier multipliers  \[
\Omega_{\la, q}(\xi)=\sum_{a\in Z_q}\sum_{\bfa\in Z^n_q}e(-\la a/q)F_q(a,\bfa) \zeta(q^2(\xi-\bfa/q))\ds(\la^{1/2}(\xi-\bfa/q))\]
where  \[
F_q(a,\bfa)=q^{-n}\sum_{s\in Z_q^n} e(|s|^2a/q+s\cdot \bfa/q).\]
is a Gaussian sum, $e(z)$ is shorthand for $e^{2\pi i z}$, $\sigma$ denotes the normalized Euclidean surface measure on the sphere of radius one in $\R^n$, and \[
\ds(\xi)=\int e(x\cdot\xi)\,d \sigma(x)\]
denotes the  $\R^n$-Fourier transform of $\sigma$. 
The function $\zeta$ is a smooth cutoff  function which is identically one on the cube $[-1/10, 1/10]^n$ and zero outside of cube $[-1/5,1/5]^n$. The notation $Z_q$ for a natural number $q$ is used for the cyclic group $\Z/q\Z$ and $U_q$ is used to denote the multiplicative group $Z_q^*$; for convenience we take $U_1=Z_1$.

Previous studies of spherical averages consider  multipliers which take the shape  \[
\sum_{q=1}^\infty \sum_{a\in U_q} m_\la^{a/q}(\xi) e(-\la a/q),\]
where $m_\la^{a/q}(\xi)$ represents the contribution of the major arc integral arising in the study of the (normalized) discrete Fourier transform\eq\label{eq00}
\frac{1}{r(\la)}\sum_{|x|^2=\la}e(x\cdot \xi)\ee 
via the circle method; the argument involving the circle method is carried out in \cite{MSW}. The characters  $e(-\la a/q)$ may appear harmless, however there is dependence on $\la$ and this in turn complicates the analysis when seeking $\ell^p$-boundedness for the maximal operator when $p<2$. Sequences which are factorial, as the one we have chosen, reduce this influence. In particular, for a given $q$ the characters are identically one for all but finitely many $\la$. As we shall see, the influence of the character can be disposed of completely. 

To proceed we notice that each  operator $A_{\la_l} f$ can be written in the form \[
\mathscr{F}^{-1}(\Omega_{l,1}\widehat{f})+\mathscr{F}^{-1}(\Omega_{l,2}\widehat{f}-\Omega_{l,1}\widehat{f})+...+\mathscr{F}^{-1}(\Omega_{l,j_0}\widehat{f}-\Omega_{l,j_0-1}\widehat{f})+(A_{\la_l}f-\mathscr{F}^{-1}(\Omega_{l,j_0}\widehat{f}))\]
where $
\Omega_{l,j}(\xi):=\Omega_{\la_l,Q_j}$ with $Q_j=2^j!.$   The notation $\widehat{f}$ is used to denote the Fourier transform on the group $\Z^n$, i.e.\[
\widehat{f}(\xi)=\sum_{x\in\Z^n}f(x)e(x\cdot\xi),\] while the notation $\mathscr{F}^{-1}$ is used for the inverse Fourier transform on the group $\Z^n$, so that \[
f(x)=\int_{\Pi^n}\widehat{f}(\xi)e(-x\cdot\xi)\,d\xi.\]

The parameter $j_0$ is chosen to be the integer $j$ such that $\la_l\in H_j$ where $H_j=[2^{2^{16(j-1)}}, 2^{2^{16j}})$.  With this selection  for $j_0$ we have  \[
\sup_{l>0} |A_{\la_l}f| \leq \sum_{i=1}^\infty \sup_{\la_l\geq 2^{2^{16(i-1)}}}  |\mathscr{F}^{-1}(\Omega_{l,i}\widehat{f}-\Omega_{l,i-1}\widehat{f})|
+\sum_{j_0=1}^{\infty}\sup_{\la_l\in H_{j_0}}|A_{\la_l}f-\mathscr{F}^{-1}(\Omega_{l,j_0}\widehat{f})|,\]
from which it immediately follows that \[
||\sup_{l>0} |A_{\la_l}f|\,||_{\ell^p} \leq \sum_{i=1}^\infty ||\sup_{\la_l\geq 2^{2^{16(i-1)}}}  |\mathscr{F}^{-1}(\Omega_{l,i}\widehat{f}-\Omega_{l,i-1}\widehat{f})|\,||_{\ell^p} 
+\sum_{j_0=1}^{\infty}||\sup_{\la_l\in H_{j_0}}|A_{\la_l}f-\mathscr{F}^{-1}(\Omega_{l,j_0}\widehat{f})|\,||_{\ell^p} .\]
It is to be understood that $\Omega_{l,0}$ is identically zero.

Theorem 1 follows once we have estimates of the form 
\eq\label{2.1}
 ||\sup_{\la_l\geq 2^{2^{16(i-1)}}}  |\mathscr{F}^{-1}(\Omega_{l,i}\widehat{f}-\Omega_{l,i-1}\widehat{f})|\,||_{\ell^p}=O(2^{-\delta i})||f||_{\ell^p} 
 \ee
 and 
 \eq\label{2.2}
\sum_{j_0=1}^{\infty}||\sup_{\la_l\in H_{j_0}}|A_{\la_l}f-\mathscr{F}^{-1}(\Omega_{l,j_0}\widehat{f})|\,||_{\ell^p} =O(2^{-\delta j_0})||f||_{\ell^p}
\ee
for each $p>1$  with positive constants\footnote {Throughout $\delta$ denotes a generic positive constant which is not necessarily the same at each appearance.} $\delta$ which are allowed to depend on $p$ and $n$. These inequalities are obtained by an interpolation between  $\ell^{1+\eps}$-estimates (for any fixed $\eps>0$) with stronger $\ell^2$ estimates. 

A discussion for each of these estimates is given in the the subsequent two sections; $\ell^{1+\eps}$-estimates are shown in the next section which is followed by a section dedicated to obtaining suitable $\ell^2$-estimates. The note then concludes with a few remarks.

\bigskip
\section{$\ell^{1+\eps}$-estimates}

We begin with an important lemma needed in the proof. The precise result we need is actually trivial in the current setting and is something that relies only on the fact that the sequence is very sparse.

\bigskip
\begin{lemma}\label{partial}
We have that \[
|| \sup_{N_0\leq \la_l\leq N_0^2}|A_{\la_l}f|\,||_{\ell^p}\leq||f||_{\ell^p}
\]
for all $p\geq1$.
\end{lemma}
\bigskip

As an immediate corollary we get that 
\[
||\sup_{\la_l\in H_{j_0}}|A_{\la_l}f|\,||_{\ell^{p}} \leq 16||f||_{\ell^{p}}\]
for each $p\geq1$. To get \eqref{2.1} it suffices for us to work with the individual multipliers $\Omega_{l, i}$.  For these multipliers we have the following lemma. 

\bigskip
\begin{lemma}\label{lem2}
When $i\geq1$ we have that\[
||\sup_{\la_l\geq 2^{2^{16(i-1)}}} |\mathscr{F}^{-1}(\Omega_{l, i} \widehat{f})|\,||_{\ell^p}\< ||f||_{\ell^p}\]
for all $p>1$. 
\end{lemma}
\bigskip

\begin{proof}
The crucial observation is that when $\la_l\geq 2^{2^{16(i-1)}}$ we have $2^{2^{2l}}\geq 2^{2^{16(i-1)}}
$, easily giving $l \geq i$ for $i>1$. Hence
\[
\Omega_{l, i}(\xi)=\sum_{a\in Z_{Q_i}}\sum_{\bfa\in Z^n_{Q_i}}F_{Q_i}(a,\bfa) \zeta(Q_i^2(\xi-\bfa/Q_i))\ds(\la_l^{1/2}(\xi-\bfa/Q_i))\]
as $e(\la_l a/Q_i)$ is identically one. And of course if $i=1$ then we clearly have $l\geq i$ so that the same statement holds. 

With the character gone, the functions  $\Omega_{l, i}(\xi)$ take the same shape as the multipliers treated in \cite{cook} and can be studied accordingly. This begins by writing  \[
\Omega_{l, i}(\xi)=\left( \sum_{\bfa\in Z^n_{Q_i}}\sum_{a\in Z_{Q_i}}F_{Q_i}(a,\bfa) \zeta(Q_i^2(\xi-\bfa/Q_i))\right)\left( \sum_{\bfa\in Z^n_{Q_i}} \zeta(Q_i(\xi-\bfa/Q_i))\ds(\la_l^{1/2}(\xi-\bfa/Q_i))\right),\]
which holds as $Q_i\geq 2$ because $1/(5Q_i^2)\leq 1/(10Q_i)$. This then gives that \[
\sup_{\la_l\geq 2^{2^{16(i-1)}}} |\mathscr{F}^{-1}(\Omega_{l, i} \widehat{f})|=\sup_{\la_l\geq 2^{2^{16(i-1)}}} |\mathscr{F}^{-1}(U_i V_{\la_l,i}\widehat{f})|\]
where \[
U_i(\xi)= \sum_{\bfa\in Z^n_{Q_i}}\sum_{a\in Z_{Q_i}}F_{Q_i}(a,\bfa) \zeta(Q_i^2(\xi-\bfa/Q_i))\]
and \[
V_{\la_l,i}(\xi)= \sum_{\bfa\in Z^n_{Q_i}} \zeta(Q_i(\xi-\bfa/Q_i))\ds(\la_l^{1/2}(\xi-\bfa/Q_i)).\]
In turn \[
\sup_{\la_l\geq 2^{2^{16(i-1)}}} |\mathscr{F}^{-1}(\Omega_{l, i} \widehat{f})| = \sup_{\la_l\geq 2^{2^{16(i-1)}}} |\mathscr{F}^{-1}(V_{\la_l,i}) * \mathscr{F}^{-1}(U_i)*f|\]
\[\leq \sup_{\la_l\geq 2^{2^{16(i-1)}}} |\mathscr{F}^{-1}(V_{\la_l,i}) * \left(\mathscr{F}^{-1}(U_i) * f\right)|.\]

Now take the $\ell^p$-norms of both sides. The maximal operators associated to $\mathscr{F}^{-1}(V_{\la_l,i}*\widehat{f})$ are bounded on $\ell^p$ for $p>1$. This follows by an application of (\cite{MSW}, Corollary 2.1) which reduces $\ell^p$-boundedness of these maximal operators to the  $L^p(\R^n)$-boundedness of the $\R^n$-spherical maximal averages over  the sequence $\la_l$ (which holds when $p>1$ as we have $n\geq 5$). Hence, when $p>1$ we have\eq\label{0001}
||\sup_{\la_l\geq 2^{2^{16(i-1)}}} |\mathscr{F}^{-1}(V_{\la_l,i}) * \left(\mathscr{F}^{-1}(U_i) * f\right)|\,||_{\ell^p}\< ||\left(\mathscr{F}^{-1}(U_i) * f\right)||_{\ell^p}.\ee
The remaining $\ell^p$ norm is treated in (\cite{cook}, Section 6), from which we have the estimate \eq\label{0002}
||\left(\mathscr{F}^{-1}(U_i) * f\right)||_{\ell^p}\< Q_{i}^{1-n}\#\{s\in Z_{Q_{i}}^n:|s|^2=0\} ||f||_{\ell^p}.\ee
when $p\geq 1$. We have  \[Q_{i}^{1-n}\#\{s\in Z_{Q_{i}}^n:|s|^2=0\} =O(1)\]
 with an implied constant independent of $i$ by (\cite{cook}, Lemma 14). Combining \eqref{0001} and \eqref{0002} completes the proof.
\end{proof}
\bigskip

The $\ell^{1+\eps}$-estimates we need to complete the outline given in section 2 are  now obtained by a direct application of the triangle inequality. 

\bigskip
\begin{lemma}\label{l^1}
For each $i\geq 1$ and  $j_0\geq 1$ we have\eq
||\sup_{\la_l\geq2^{2^{16(i-1)}}}|\mathscr{F}((\Omega_{l, i}-\Omega_{l, i-1})\widehat{f})|\,||_{\ell^p}\<   ||f||_{\ell^p}\ee
and 
\eq
||\sup_{\la_l\in H_{j_0}}|(A_{\la_l}f-\mathscr{F}(\Omega_{l, j_0}\widehat{f}))|\,||_{\ell^p}\< ||f||_{\ell^p}\ee 
 for all $p>1$. The implied constants do not depend on $i$ or $j_0$. 
\end{lemma}

\bigskip


\section{$\ell^2$-estimates}

Our goal in this section is to show a suitable result which, as previously mentioned, when interpolated against Lemma \ref{l^1} implies Theorem \ref{thm1}. The statement we need is as follows.

\bigskip
\begin{lemma}
The inequality 
\[
||\sup_{\la_l\geq2^{16^{i-1}}}|\mathscr{F}((\Omega_{l, i}-\Omega_{l, i-1})\widehat{f})|\,||_{\ell^2}\<  2^{-\delta i} ||f||_{\ell^2}\]
and the inequality
\[
||\sup_{\la_l\in H_{j_0}}|(A_{\la_l}f-\mathscr{F}(\Omega_{l, j_0}\widehat{f}))|\,||_{\ell^2}\< 2^{-\delta j_0} ||f||_{\ell^2}.\] 
hold for each $i\geq 1$ and  $j_0\geq 1$  with implied constants independent of $i$ and $j_0$.
\end{lemma}
\bigskip

Estimates at $\ell^2$ are obtained through  the introduction of the multipliers \[
M_{\la,h}(\xi)=\sum_{q\in I_h}\sum_{a\in U_q}\sum_{\bfa\in Z_\bfq^n}e(-\la a/q) F_q(a,\bfa)\zeta(10^h(\xi-\bfa/q))\ds(\la^{1/2}((\xi-\bfa/q)) \]
where $I_h=[2^h,2^{h+1})$ for $h\geq0$. These terms are comparable to (dyadic portions of) the multipliers that appear in the more standard approximations for the normalized Fourier transforms \[
\widehat{\omega_\la}(\xi):=\frac{1}{r(\la)}\sum_{|x|^2=\la}e(x\cdot\xi)\]
that are employed in \cite{MSW} and subsequent works, the major distinction here is that the cutoff $\zeta$  has been made uniform on each dyadic scale. The precise statement of our approximation to $\widehat{\omega_\la}$ using these terms is given as the next result. This lemma is in fact identical to (\cite{cook}, Lemma 6), however we shall include a proof. 

\bigskip
\begin{lemma} 
For every $\la\geq1$ we have \[
\widehat{\omega_\la}(\xi)=\sum_{h=0}^{\infty} M_{\la,h}(\xi)+O(\la^{-\delta})\]
uniformly in $\xi$.
\end{lemma}
\bigskip
\begin{proof}

For fixed $h\geq0$ and $q\in I_h$ we look at  \[
\sum_{a\in U_q}\sum_{\bfa\in Z_q^n}e(-\la a/q)F_q(a,\bfa)\left(\zeta(q(\xi-\bfa/q)-\zeta(10^h(\xi-\bfa/q))\right)\ds(\la^{1/2}(\xi-\bfa/q)).\] 
When $\xi$ is fixed there is at most one $\bfa$ for which $\zeta(q(\xi-\bfa/q)-\zeta(10^h(\xi-\bfa/q))$ is nonzero, and looking on the support of any of these terms we see that $|\xi-\bfa/q|\geq 10^{-h}.$ We have well known Fourier decay estimates for $\widetilde{d\sigma}$  that are useful when $h$ is small in terms of $\la$. More precisely, if  $h<\delta \log\,\la$  we  have  that 
\[
|\sum_{a\in U_q}\sum_{\bfa\in Z_q^n}e(-\la a/q)F_q(a,\bfa)\left(\zeta(q(\xi-\bfa/q)-\zeta(10^l(\xi-\bfa/q))\right)\ds(\la^{1/2}(\xi-\bfa/q))|\]\[\< \left(\sum_{a\in U_q}\sup_{\bfa\in Z_q}|F_q(a,\bfa)|\right)\la^{-\delta}\]\[
\< q^{1-n/2}\la^{\delta}\]
holds uniformly in $\xi$ thanks to any Fourier decay estimate of the form \[
|\ds(\la^{1/2}(\xi-\bfa/q))|\< (\la/10^h)^{-c}\]
where the precise value of $c$ is not currently relevant (although the correct value $c=(n-1)/2$ will be used below).  
Now we may sum in  $q\geq1$ when $n>4$ and arrive at the  upper bound \[
\sum_{h<\delta\log\,\la}\sum_{q\in I_h}\sum_{a\in U_q}\sum_{\bfa\in Z_q^n}e(-\la a/q)F_q(a,\bfa)\left(\zeta(q(\xi-\bfa/q)-\zeta(10^h(\xi-\bfa/q))\right)\ds(\la^{1/2}(\xi-\bfa/q))\]\[=O(\la^{-\delta}).\] 

We must also consider the remaining terms with $h\geq\delta \log\,\la$ with $\delta$ the constant used directly above, specifically this is for 
\[
\sum_{h\geq\delta\log\,N}\sum_{q\in I_h}\sum_{a\in U_q}\sum_{\bfa\in Z_q^n}e(-\la a/q)F_q(a,\bfa)\left(\zeta(q(\xi-\bfa/q)-\zeta(10^h(\xi-\bfa/q))\right)\ds(\la^{1/2}(\xi-\bfa/q)).\] 
Again, appealing to the disjointness of the supports of the summands over the various $\bfa\in Z_q$ we see the upper bound  \[
\sum_{\bfa\in\Z^n}e(-\la a/q)F_q(a,\bfa)\left(\zeta(q(\xi-\bfa/q)-\zeta(10^h(\xi-\bfa/q))\right)\ds(\la^{1/2}(\xi-\bfa/q))\< \sup_{\bfa\in Z_q} |F_q(a,\bfa)|\]
uniformly in $\xi\in\Pi^n$ because\[
\left(\zeta(q(\xi-\bfa/q)-\zeta(10^h(\xi-\bfa/q))\right)\ds(\la^{1/2}(\xi-\bfa/q))\<1.\] Using the Gauss sum estimate  $|F_q(a,\bfa)|\<q^{-n/2}$ reduces the consideration to \[
\sum_{h\geq\delta\log\,N}\sum_{q\in I_h}\sum_{a\in U_q}q^{-n/2}\leq\sum_{q\geq \la^{\delta}}q^{1-n/2}.\]
This is $O(\la^{-\delta})$ as $n>4$ and the proof is complete. 
\end{proof}
\bigskip

We consider approximations for the $\Omega_{l, j}(\xi)$ terms by the $M_{\la_l,h}(\xi)$. The precise form is taken from \cite{cook} and takes the shape \eq\label{3.1}
\Omega_{l, j}(\xi)=\sum_{h=0}^{j-1}M_{\la_l,h}(\xi)+E^{(1)}_{l,j}(\xi)+E^{(2)}_{l,j}(\xi)\ee
where \[
E^{(1)}_{l,j}(\xi)=\sum_{h=0}^{j-1}\sum_{q\in I_h}\sum_{a\in U_q}\sum_{\bfa\in Z^n_q}e(-a\la_l/q)F_q(a,\bfa)\left(\zeta(Q_j^2(\xi-\bfa/q))-\zeta(10^h(\xi-\bfa/q))\right)\ds(\la_l^{1/2}(\xi-\bfa/q))\]
and \[
E^{(2)}_{l,j}(\xi)=\sum_{q|Q_j;\,q\geq 2^j}\sum_{a\in U_q}\sum_{\bfa\in Z^n_q}e(-a\la_l/q)F_q(a,\bfa)\zeta(Q_j^2(\xi-\bfa/q))\ds(\la_l^{1/2}(\xi-\bfa/q)).\]
In turn we have that \[
\Omega_{l, j}(\xi)-\Omega_{l, j-1}(\xi)=M_{\la_l,j-1}(\xi)+E^{(1)}_{l,j}(\xi)+E^{(2)}_{l,j}(\xi)-E^{(1)}_{l,j-1}(\xi)-E^{(2)}_{l,j-1}(\xi).\]


The derivation of \eqref{3.1} makes repeated use of the basic transformation \[
\sum_{a\in\Z_q}g(a/q)\leftrightarrow \sum_{d|q} \sum_{a\in U_q}g(a/d)\]
and a relatively easy observation about exponential sums (\cite{cook}, Lemma 2); a more detailed discussion can be found in (\cite{cook}, section 2.4).

We now present the necessary $\ell^2$-estimates that complete the proof of our main result.

\bigskip
\begin{lemma}
For each $l\geq1$ and $i\geq1$ the following estimates hold:\[
 ||\sup_{\la_l\geq 2^{2^{16(i-1)}}}  |\mathscr{F}^{-1}(E^{(1)}_{l,i}\widehat{f})|\,||_{\ell^2}=O(2^{-\delta i}); \]
 \[
 ||\sup_{\la_l\geq 2^{2^{16(i-1)}}}  |\mathscr{F}^{-1}(E^{(2)}_{l,i}\widehat{f})|\,||_{\ell^2}=O(2^{-\delta i}); \]
and  \[
 ||\sup_{\la_l\geq 2^{2^{16(j-1)}}}  |\mathscr{F}^{-1}(M_{\la_l,j}|\,||_{\ell^2}=O(2^{-\delta j}). \]
\end{lemma}
\bigskip

\begin{proof}
The inequalities are treated in order of appearance. We have  \[
 ||\sup_{\la_l\geq 2^{2^{16(i-1)}}}  |\mathscr{F}^{-1}(E^{(1)}_{l,i}\widehat{f})|\,||_{\ell^2}\]\[
 \leq \sum_{\la_l\geq 2^{2^{16(i-1)}}} || \mathscr{F}^{-1}(E^{(1)}_{l,i}\widehat{f})||_{\ell^2}\]\[
 \leq \sum_{\la_l\geq 2^{2^{16(i-1)}}} \sum_{h=0}^{i-1}\sum_{q\in I_h}\sum_{a\in U_q} || \mathscr{F}^{-1}(\sum_{\bfa\in Z^n_q}F_q(a,\bfa)\left(\zeta(Q_i^2(\cdot-\bfa/q))-\zeta(10^h(\cdot-\bfa/q))\right)\ds(\la^{1/2}(\cdot-\bfa/q))\widehat{f})||_{\ell^2}\]\[
 = \sum_{\la_l\geq 2^{2^{16(i-1)}}} \sum_{h=0}^{i-1}\sum_{q\in I_h}\sum_{a\in U_q} || \sum_{\bfa\in Z^n_q}F_q(a,\bfa)\left(\zeta(Q_i^2(\cdot-\bfa/q))-\zeta(10^h(\cdot-\bfa/q))\right)\ds(\la^{1/2}(\cdot-\bfa/q))\widehat{f}||_{L^2(\Pi^n)}\]
When we consider \[
\sum_{\bfa\in Z^n_q}F_q(a,\bfa)\left(\zeta(Q_i^2(\cdot-\bfa/q))-\zeta(10^h(\cdot-\bfa/q))\right)\ds(\la^{1/2}(\cdot-\bfa/q))\]
we see that, with fixed $\xi$, $a$,  and $q$,  there is at most one term which is nonzero as the functions  \[
\zeta(Q_i^2(\xi-\bfa/q))-\zeta(10^l(\xi-\bfa/q)\]
 have disjoint supports as $\bfa$ varies.  Thus we now need to show that\[
 \sum_{\la_l\geq 2^{2^{16(i-1)}}}\sum_{h=0}^{i-1}\sum_{q\in I_h}\sum_{a\in U_q}\sup_{\bfa\in Z^n_q}\sup_{\xi\in\Pi^n}|F_q(a,\bfa)\left(\zeta(Q_i^2(\xi-\bfa/q))-\zeta(10^h(\xi-\bfa/q))\right)\ds(\la^{1/2}(\xi-\bfa/q))|\]\[
\leq \sum_{\la_l\geq 2^{2^{16(i-1)}}}\sum_{h=0}^{i-1}\sum_{q\in I_h}\sum_{a\in U_q}q^{-n/2}\sup_{\bfa\in Z^n_q}\sup_{\xi\in\Pi^n}|\left(\zeta(Q_i^2(\xi-\bfa/q))-\zeta(10^h(\xi-\bfa/q))\right)\ds(\la^{1/2}(\xi-\bfa/q))|\]\[
=O(2^{-\delta i}).\]

 We now apply the Fourier decay estimate for $\ds$  as before (and now with the precise constant $c=-(n-1)/2$), giving  \[
\sum_{\la_l\geq 2^{2^{16(i-1)}}}\sum_{q=1}^{2^i}\sum_{a\in U_q}q^{-n/2}
\left(\frac{1}{1+\la_l^{1/2}/Q^2_i}\right)^{(n-1)/2}\]
\[
\<\sum_{\la_l\geq 2^{2^{16(i-1)}}}
\left(\frac{Q^2_i}{\la_l^{1/2}}\right)^{(n-1)/2}\]\[
=\sum_{\la_l\geq 2^{2^{16(i-1)}}}
\left(\frac{(2^i!)^4}{(2^l!)}\right)^{(n-1)/4}\]\[
\leq\sum_{l\geq 8(i-1)}\left(\frac{(2^i!)^4}{(2^l!)}\right)^{(n-1)/4}\]\[
=\sum_{l\geq 8(i-1)}\left(\frac{(2^i!)^4}{(2^{4i}!)}\frac{1}{(2^{4i}+1)(2^{4i}+2)...2^{l}}\right)^{(n-1)/4}\]\[
\leq\sum_{l\geq 8(i-1)}\left(\frac{1}{(2^{4i})^{l-4i}}\right)^{(n-1)/4}\]
\[\leq\sum_{l\geq 4i-8}\left(\frac{1}{(2^{4i})^{l}}\right)^{(n-1)/4}\]
\[\leq\sum_{l\geq 4i-8}2^{-i(n-1)l}\]
\[=O((2^{-\delta i})\]
as long as $i\geq3$. The $i=1,2$ cases follow by adjusting the implied constant as necessary, and then the inequality as stated holds.

We now consider the inequality for the  $E^{(2)}_{l,i}$ terms. Here we have \[
 ||\sup_{\la_l\geq 2^{2^{16(i-1)}}}  |\mathscr{F}^{-1}(E^{(2)}_{l,i}\widehat{f})|\,||_{\ell^2}\]\[= ||\sup_{\la_l\geq 2^{2^{16(i-1)}}}  |\mathscr{F}^{-1}(\sum_{q|Q_i;\,q\geq 2^j}\sum_{a\in U_q}\sum_{\bfa\in Z^n_q}e(-a\la_l/q)F_q(a,\bfa)\zeta(Q_i^2(\cdot-\bfa/q))\ds(\la^{1/2}(\cdot-\bfa/q))\widehat{f})|\,||_{\ell^2} \]\[
 \leq 
 \sum_{q|Q_i;\,q\geq 2^i}\sum_{a\in U_q}||\sup_{\la_l\geq 2^{2^{16(i-1)}}}  |\mathscr{F}^{-1}(\sum_{\bfa\in Z^n_q}F_q(a,\bfa)\zeta(Q_i^2(\cdot-\bfa/q))\ds(\la^{1/2}(\cdot-\bfa/q))\widehat{f})|\,||_{\ell^2}.\]
 The right hand side of this inequality is bounded by a multiple of\[
 \sum_{q|Q_i;\,q\geq 2^i}\sum_{a\in U_q}||  \mathscr{F}^{-1}(\sum_{\bfa\in Z^n_q}F_q(a,\bfa)\zeta(Q_i^2(\cdot-\bfa/q))\widehat{f})||_{\ell^2},\]
 which follows by applying  (\cite{MSW}, Corollary 2.1) as was done in Lemma \ref{lem2}. We also have \[
 ||  \mathscr{F}^{-1}(\sum_{\bfa\in Z^n_q}F_q(a,\bfa)\zeta(Q_i^2(\cdot-\bfa/q))\widehat{f})||_{\ell^2}\<\sup_{\bfa\in Z_q}|F_q(a,\bfa)|\,||f||_{\ell^2}\]
 by Plancherel's equality and the disjoint supports of the $\zeta(Q_i^2(\cdot-\bfa/q))$ terms. The inequality we claimed above follows as \[
 \sum_{q|Q_i;\,q\geq 2^i}\sum_{a\in U_q}q^{-n/2}\]\[
 \leq \sum_{q\geq 2^i}q^{1-n/2}\]\[
 =O(2^{-\delta i})\]
 when $n>4$.
 
 The remaining inequality  is established in a similar manner as is in the argument used above for deducing the inequality for the   $E^{(2)}_{l,i}$ terms. Indeed, the key observation here is again that the exponential sums provide necessary gains as the $q$ are restricted to the regime $q\geq2^{j-1}$.
\end{proof}

\bigskip

\section{Further remarks}

The discussion here follows along the lines of the argument in \cite{cook} and a byproduct of this is that, in contrast to the arguments of \cite{Kevin} and \cite{Lacey},  there is no reliance on estimates for Kloostermann sums. As a result one can easily extend our result to more general positive definite integral forms of higher degree like those considered in \cite{Akos}. The statements and proofs require only a minor modification of what is presented below. The sequence  $\la_l$ needs to modified  due to the fact that the  polynomials in question may not be $universal$ in the sense that they may not represent all large or even all sufficiently large $\la\in\mathbb{N}$. A consequence of the results  in \cite{Bi} is that a nonsingular positive definite form, say $\p$,  of degree $d$ in at least $(d-1)2^d$ variables represents all $\la$ in a given infinite arithmetic progression $\Gamma=\{al+b:l\in\Z, x\geq l_0\}$; details can be found in (\cite{Akos}, Lemma 10). 

\bigskip
\begin{thm}\label{thm2}
Let $\p$ be a positive definite integral form of degree $d>1$ with $\mathcal{B}(\p)\geq(d-1)2^d$ and let $\Gamma$ (as above) be a sequence of regular values of $\p$. Define $\la_l= 2^l!+b$. The maximal operator \[
\sup_{l\geq l_0}\frac{1}{\la_l^{n/d-1}}\left|\sum_{\p(x)=\la_l} f(y-x)\right|\] 
is bounded operator on $\ell^p$ for all $p>1$.
\end{thm}
\bigskip

The proof attached to such  statements requires only a very minor modification to the argument given above. Namely, at each appearance of quantities of the type $2^{2^{16i}}$ one needs to substitute something of the form $2^{2^{c_d i}}$ for a suitable constant $c_d$ which depends only on the degree of the relevant polynomial. The reader that has followed the argument above should have no trouble verifying this.

An interesting question in a slightly different direction would be to consider Theorem 1 for suitable sequences which are less sparse, namely in cases where Lemma \ref{partial} is no longer trivial. For example, it seems quite plausible that the result of Theorem \ref{thm1} should hold with the sequence $\la_l=l!$  in place of $\la_l=2^l!$.


\vskip0.2in
\noindent \author{\textsc{Brian Cook}}\\
Department of Mathematics \\
Kent State University\\
Kent, OH, USA\\
Electronic address: \texttt{briancookmath@gmail.com}


\end{document}